\documentclass[a4paper,leqno,12pt,twoside]{amsart}

\usepackage{amsmath,amscd,amssymb,upgreek} 
\usepackage[all]{xy}
\usepackage{graphicx}
\newtheorem{Lem}{Lemma}
\newtheorem{Th}[Lem]{Theorem }
\newtheorem{Cor}[Lem]{Corollary }
\newtheorem{Def}[Lem]{Definition }

\newtheorem{Prop}[Lem]{Proposition }
\newtheorem{Rem}[Lem]{Remark }
\newtheorem{Lem and Def}[Lem]{Lemma and Definition \thesection.\!\!}

\newtheorem{Prop and Def}[Lem]{Proposition and Definition   \thesection.\!\!}

 \def\cal{\mathcal}

\newcommand{\C}{\mathbb{C}}
\newcommand{\N}{\mathbb{N}}
\newcommand{\R}{\mathbb{R}}

\newcommand{\Z}{\mathbb{Z}}

\def\Supp{\mathop{\rm Supp}}
\def\log{\mathop{\rm log}}
\def\trace{\mathop{\rm trace}}

\begin{document}
\title{On the K\"ahler Rank of Compact Complex Surfaces}
\author{ Matei TOMA}  
\date{\today}

\thanks{ I wish to thank  the
Max-Planck-Institut f\"ur Mathematik in Bonn and
the University of Osnabr\"uck for their hospitality and for financial support during the preparation of this paper. Furthemore I  thank Charles Favre and Karl Oeljeklaus for useful discussions and 
the referee for his remarks  which improved the presentation of the paper.  \\
Keywords: compact complex surface, global spherical shell, closed positive current, iteration of polynomial maps.\\
Mots-cl\'e: surface complexe compacte, coquille sph\'erique globale, courant positif ferm\'e, it\'eration des applications polyn\^omiales\\
AMS Classification (2000): 32J15, 32H50 }

\address{
Institut de Math\'ematiques Elie Cartan,
Nancy-Universit\'e,
B.P. 239,
54506 Vandoeuvre-l\`es-Nancy Cedex,
France  and 
Institute of Mathematics of the Romanian Academy.}
\email{toma@iecn.u-nancy.fr}
\urladdr{http://www.iecn.u-nancy.fr/$\sim$toma/}

\begin{abstract}
 Harvey and Lawson introduced the K\"ahler rank and computed it in connection to the cone of positive exact currents
of bidimension $(1,1)$ for many classes of compact complex surfaces.
   In this paper  we extend these computations to the only further known class of surfaces not considered by them,
 that of Kato surfaces.  Our main tool is the reduction to the dynamics of associated  holomorphic contractions
 $(\C^2,0)\to (\C^2,0)$.\\
 
 Harvey et Lawson ont introduit et calcul\'e le rang de K\"ahler 
 en relation avec le c\^one des courants positifs ferm\'es  
 de bidimension $(1,1)$ pour beaucoup de classes de surfaces complexes compactes.
 Dans ce travail nous \'etendons ces calculs \`a la seule classe de surfaces connues
 et qui n'avait pas \'et\'e consid\'er\'ee par eux, celle des surfaces de Kato. 
 Notre outil principal est la r\'eduction \`a la dynamique des contractions holomorphes 
 $(\C^2,0)\to (\C^2,0)$ associ\'ees.

 \end{abstract}

\setcounter{section}{0}
\noindent
\maketitle

\section{Introduction}

In  \cite{HL83} Harvey and Lawson give a characterisation of K\"ahlerianity for compact complex surfaces
in terms of existence (or rather non-existence) of closed positive currents which are $(1,1)$-components of a boundary. The authors also investigate and describe the cones formed by such currents for many types of non-K\"ahler
surfaces: elliptic, Hopf, Inoue. Later Lamari proved that every  non-K\"ahler
surface admits non-trivial positive $d$-exact currents of bidimension $(1,1)$; cf. \cite{Lam99}.
In order to estimate the degree of non-K\"ahlerianity of a compact complex surface 
smooth positive $d$-exact currents 
are considered in  \cite{HL83} and the K\"ahler rank 
is defined as follows: the K\"ahler rank  is 
{\em two}  if the surface admits some K\"ahler metric, {\em one} if it admits some positive $d$-exact $(1,1)$-form with some supplementary property and {\em zero} in the remaining case; see the more precise Definition \ref{defKr}.
Unfortunately it is not clear whether the K\"ahler rank is a bimeromorphic invariant.

We consider in this paper instead a bimeromorphic invariant which we call the {\em modified K\"ahler rank} and which we define to be {\em two}  if the surface is  K\"ahler, {\em one} if the cone of positive $d$-exact currents of bidimension $(1,1)$ is larger than a half-line and {\em zero} if this cone is a half-line. The two notions agree in the cases considered in \cite{HL83}. Our main result is the computation of the cones of positive exact 
$(1,1)$-currents for Kato surfaces. These are surfaces whose minimal models have positive second Betti number and admit a global spherical shell (see Definition \ref{gss}) and are the only ``known'' compact complex surfaces not considered in \cite{HL83}. 
(We refer the reader to \cite{BHPV} for the general theory of compact complex surfaces.)
It  turns out that the modified  K\"ahler rank does not coincide with the K\"ahler rank
in general. In order to perform our computations we reduce ourselves to the investigation of plurisubharmonic functions with a certain invariance property with respect to polynomial automorphisms of $\C^2$ associated to 
 Kato surfaces.
 As a corollary we obtain
 
 \begin{Th}\label{foliated}
 (a) Every positive $d$-exact $(1,1)$-current on a Kato surface, and more generally on any
 ``known'' non-K\"ahlerian compact complex surface, is a foliated current for some holomorphic foliation of the surface.

 (b) All positive $d$-closed $(1,1)$-currents on  the ``known'' non-K\"ahlerian compact complex surfaces excepting on parabolic Inoue surfaces are foliated currents for some holomorphic foliations.
\end{Th} 

See Section \ref{sectKato} for definitions.

\section{The K\"ahler rank}

We let $X$ be a compact complex surface and denote by $P_{bdy}=P_{bdy}(X)$, $P_{bdy}^{\infty}=P_{bdy}^{\infty}(X)$
the cones of positive currents of bidimension $(1,1)$ which are boundaries (i.e. are $d$-exact), respectively smooth such currents. These objects have been first studied in \cite{HL83}. In \cite{Lam99} it was shown that $X$ is K\"ahler if and only if $P_{bdy}(X)$ is trivial. It is easy to see that on a non-K\"ahler surface every positive $d$-closed differential form of type $(1,1)$ is $d$-exact, cf. \cite{Lam99}. The following definition of Harvey and Lawson gives a measure of non-K\"ahlerianity by looking at the positive closed differential$(1,1)$-forms on $X$.
It is roughly speaking the largest generic rank of a positive closed $(1,1)$-form on $X$.

\begin{Def} \label{defKr}
Let $B(X)=\{x\in X \ | \ \exists\phi\in P_{bdy}^{\infty}(X) \ \phi_x\ne 0\}$.
The {\em  K\"ahler rank} of $X$ is defined to be {\em two}  if $X$ admits a K\"ahler metric. If $X$ is non-K\"ahler
the K\"ahler rank is set to be {\em one} when $B(X)$ contains a non-trivial Zariski open subset of $X$ and  {\em zero}
otherwise.
\end{Def}
It is not known whether the K\"ahler rank is a bimeromorphic invariant.

We propose also the following:
\begin{Def} 
The {\em  modified K\"ahler rank} of $X$ is defined to be {\em two}  when $X$ admits no non-trivial positive exact current of bidimension  $(1,1)$, {\em zero} when it admits exactly one such current up to a multiplicative constant and  {\em one}
otherwise.
\end{Def}

One can show easily that the modified K\"ahler rank is a bimeromorphic invariant by taking push-down and pull-back of currents through blowing up maps. 
See the proof of Proposition \ref{bimero} for a more precise description.

For elliptic surfaces, primary Hopf surfaces and Inoue surfaces  one sees that
the  K\"ahler rank and the modified K\"ahler rank coincide 
using the precise description of $P_{bdy}$ given in \cite{HL83} in these cases.

\begin{Prop}\label{prop}
 Let $T$ be a positive exact current of bidimension  $(1,1)$ on the compact complex surface $X$.
Then there is a representation $\rho:\pi_1(X)\to(\R,+)$ and a plurisubharmonic function $u$
on the universal cover
$\tilde X$ of $X$ such that $T=dd^cu$ and $u\circ g=u+ \rho(g)$ for all $g\in \pi_1(X)$. 

The function $u$ can be chosen to be smooth if $T$ is smooth.
\end{Prop} 
\begin{proof}
One has $b_1=\dim_{\C}H^1(X,\C)=\dim_{\C}H^1(X,\cal{O}_X)+\dim_{\C}H^0(X,\Omega^1_X)=h^{0,1}+h^{1,0}$, cf \cite{BHPV}.
Denoting the sheaf of closed differential $(1,0)$-forms 
by $d\cal{O}_X$ and looking at the long exact cohomology sequence of
$$0\to\C_X\to\cal{O}_X\to d\cal{O}_X\to 0$$
one gets an exact sequence
$$0\to H^0(d\cal{O}_X)\to H^1(\C_X)\to H^1(\cal{O}_X).$$
By the above equality follows now the surjectivity of the natural map $H^1(\C_X)\to H^1(\cal{O}_X).$ 
This is given by mapping a de Rham cohomology class $[\beta]$ 
of a differential form $\beta=\beta^{1,0}+\beta^{0,1}$ onto the 
Dolbeault cohomology class of $\beta^{0,1}$. 

Let now $T$ be a positive exact current of bidimension $(1,1)$ on $X$. Then $T=dS$ with
$S=S^{1,0}+S^{0,1}$,
$S^{1,0}$, $S^{0,1}$ currents of order zero and bidegree $(1,0)$ and $(0,1)$ respectively, 
and $S^{1,0}=\bar S^{0,1}$.
Since $\bar\partial S^{0,1}=0$, $S^{0,1}$ represents a cohomology class in  $H^1(\cal{O}_X)$ and let 
$\beta=\beta^{1,0}+\beta^{0,1}$ be a closed differential form with $[\beta^{0,1}]=[S^{0,1}]$ in $H^1(\cal{O}_X)$ and $U$ a current of degree $0$ on $X$ with $S^{0,1}=\beta^{0,1}+\bar\partial U$.
The lift $\tilde\beta$ of $\beta$ to the universal cover $\tilde X$ is $d$-exact and 
let $f$ be a smooth function on 
$\tilde X$ with $df=\beta$. In particular $\bar\partial f=\beta^{0,1}$.
This implies $S^{0,1}=\bar\partial (f+U)$ and 
$T=dS=d(\bar\partial (f+U)+\partial(\bar f+\bar U))=i\partial\bar\partial(2\Im m(f+U))$.

Moreover for $g\in\pi_1(X)$ we have $d(f\circ g-f)=0$ hence $f\circ g-f$ must be constant.
Set $\rho(g)=2\Im m(f\circ g-f)$.
The current $2\Im m(f+U)$ is associated to a plurisubharmonic function $u$ on $\tilde X$. Since 
$u-2\Im mf=2\Im mU$ comes from $X$ we see that $ u$ has the desired
 automorphy behaviour with respect to the action of 
 $\pi_1(X)$.

It is clear that $u $ can be chosen to be smooth when $T$ is smooth.
\end{proof}

\begin{Def} 
We say that an effective reduced divisor $C=C_1+...+C_n$ on $X$ is a {\em cycle of rational curves}
if $n\ge1$, $C_1$,..., $C_n$ are rational curves and either $n=1$ and $C_1$ has a node or
$n>1$, all components  $C_1$,..., $C_n$ are smooth and the dual graph of $C$ is cyclic.
\end{Def}

\begin{Cor}
 For a compact complex surface $X$ with a cycle of rational curves
$C$ and $b_1(X)=1$ the K\"ahler rank is zero.
\end{Cor} 
\begin{proof}
Under the above hypotheses the natural map $\Z\cong\pi_1(C)\to\pi_1(X)$ is an isomorphism by a Theorem of Nakamura, 
\cite{Nak90}.
Let $g$ be a generator of $\pi_1(X)$. If the K\"ahler rank of $X$  were one,
we would get a smooth non-constant plurisubharmonic function $u$ on the universal cover $\tilde X$ satisfying
$$u\circ g=u+c$$
for some constant $c\in\R$ by Proposition \ref{prop}.
The inverse image $\tilde C$ of the cycle of rational curves $C$ of $X$ is an infinite chain of rational curves
on   $\tilde X$. Since 
$u$ is smooth and constant on each link of this chain we would get $c=0$ hence  $ u$ would descend to $X$.
But here  $u $ must be constant contradicting our assumptions. 
\end{proof}

\begin{Prop}\label{bimero}
 The modified K\"ahler rank is a bimeromorphic invariant.
\end{Prop} 
\begin{proof}
It is known that for smooth surfaces the property of being K\"ahler is invariant under bimeromorphic
transformations. 

Let now  $X$ be non-K\"ahler and  $X'$ the blown-up surface at a point $x\in X$.
Any plurisubharmonic function $u$ on the universal cover $\tilde X$ may be lifted to $\tilde X'$.
Conversely let  $ u'$  be a plurisubharmonic function on $\tilde X'$. It will be constant on the exceptional divisors of $\tilde X'$ coming from $X'\to X$. Thus $ u'$  is the pull-back of a plurisubharmonic function on $\tilde X$, which satisfies the same invariance condition as $u$ with respect to the action of $\pi_1(X)\cong \pi_1(X')$. 
This gives a bijective correspondence between $P_{bdy}(X)$ and $P_{bdy}(X')$, hence the invariance of the
modified K\"ahler rank. 
\end{proof}

From the above argument it also follows that a counterexample to the bimeromorphic invariance of the K\"ahler rank could only be given by some non-elliptic non-K\"ahler minimal surface $X$ of K\"ahler rank zero admitting
 a continuous plurisubharmonic function $u$ on the universal cover $\tilde X$ which
is smooth outside a discrete subset of 
$\tilde X$ and  exhibits the automorphy behaviour from Proposition \ref{prop}. 
From \cite{HL83}
we gather that this is not the case of Hopf surfaces or Inoue surfaces and as we shall see it is not going to be the case of Kato surfaces either.

Finally let us mention the following description of the cone $P(X)$ of positive $d$-closed $(1,1)$-currents:
\begin{Rem}\label{rem}
Let $C_1$, ..., $C_n$ be the irreducible curves of negative self-intersection on a
 non-K\"ahlerian surface $X$ and $[C_1]$, ..., $[C_n]$ the corresponding currents of integration. Then
 $$P(X)=P_{bdy}(X)+\sum_i\R_{\ge 0}[C_i].$$
\end{Rem}
\begin{proof}
Let us denote by $F_j$ the irreducible curves of $X$ of zero self-intersection. Since $X$ is non-algebraic it will admit no curve of positive self-intersection. Thus the Siu decomposition of a positive closed current $T$ takes the form:
$$T=\sum_i a_i[C_i]+\sum_j b_j[F_j]+ R,$$
where $a_i,b_j\in\R$ and $R$ is a positive current whose Lelong level sets are finite.
By  \cite{Lam99'} Prop. 4.3 and its proof   $R$ is nef. Since $X$ is non-K\"ahlerian it follows that $R$ must be exact, cf. \cite{Lam99'} Thm. 7.1.  The integration currents
$[F_j]$ are also exact, since $F_j^2=0$ and the intersection form on $H^{1,1}_{\R}(X)$ is negative definite for a non-K\"ahlerian surface.
\end{proof}


\section{Kato surfaces}\label{sectKato}

We recall here some facts on Kato surfaces which will be needed in the sequel. We refer the reader to 
\cite{Kato78}, \cite{Dl84}, \cite{Dl88}, \cite{DO99}, \cite{Fav00} for more details.


\begin{Def} \label{gss}
A {\em global spherical shell} on a compact complex surface $X$ is the image $\Sigma$ of the sphere $S^3$ by a holomorphic embedding of a neighbourhood of  $S^3$ from $\C^2$ into $X$ such that $X\setminus\Sigma$ is connected. 
We call $X$ a Kato surface if $X$ admits a global spherical shell 
and the second Betti number of a minimal model of $X$ is positive.
\end{Def}

The other minimal surfaces admitting a global spherical shell are the Hopf surfaces. 
Their second Betti number is zero.

The notion of global spherical shell and a construction method for surfaces with 
global spherical shells were introduced by Ma. Kato in \cite{Kato78}. 
An important analytic object associated to this 
construction method was introduced and studied by Dloussky in  \cite{Dl84}. It is the germ of a holomorphic map
$(\C^2,0)\to (\C^2,0)$. We shortly recall the facts now.

Take the unit ball $B$ in $\C^2$ around the origin and blow-up the origin.
Choose a point $P_1$ on the exceptional curve $C_1$ thus obtained and blow it up again.
Continue by blowing up a point on the last created exceptional curve. After $n$ blow-ups one considers
the blowing down map $\pi:B'\to B$. 
The exceptional divisor on $B'$ is a tree of $n$ smooth rational curves. The only 
$(-1)$-curve among them is the last created curve $C_n$. 
Choose a point $P_n$ on $C_n$, a biholomorphic map
$\sigma:\bar B\to\sigma(\bar B)$ onto a small compact neighbourhood of $P_n$ in $B'$ and glue the two
components of the boundary of 
$B'\setminus\sigma(B)$ by means of $\sigma\circ\pi$.
In this way one obtains a minimal compact complex surface $X$.
One can show that the image of $S^3$ through $\sigma$ is a global spherical shell on $X$
and that $b_2(X)=n$. Thus $X$ is a minimal Kato surface. 

The images of the exceptional curves are the only (compact) rational curves on $X$. They form an effective reduced
divisor which we denote by $D$.
Depending on the structure of $D$ one subdivides the class of minimal Kato surfaces 
into:
\begin{enumerate}
\item {\it Enoki surfaces}, 
when $D$ is a cycle of rational curves and $D$ is homologically trivial,
\item {\it intermediate surfaces}, when $D$ consists of
a cycle of rational curves and of at least
one further rational curve attached to the cycle,
\item {\it Inoue-Hirzebruch surfaces}, when $D$ consists of
one or two cycles of rational curves and $D$ is not homologically trivial.
\end{enumerate}
In particular Kato surfaces admit cycles of rational curves and therefore their K\"ahler rank is zero.

In the case of Enoki surfaces a further curve might appear. In such a case the curve will be elliptic and the surface is called a {\it parabolic Inoue surface}.
 
By definition the {\em Dloussky germ} associated to $X$ is the germ of the map 
$f:=\pi\circ\sigma:B\to B$ around the origin. It can be shown that the conjugacy class of this germ determines 
the isomorphy class of $X$. We shall denote by $X(f)$ a Kato surface associated to such a germ of holomorphic map
$f:(\C^2,0)\to (\C^2,0)$. 
One can relate certain analytic objects on $X$ to germs of objects on $(\C^2,0)$ which are invariant under
$f$ as follows.

One recovers first the universal cover $\tilde X$ of $X$ from the above construction 
by considering an infinite number of copies $(A_i)_{i\in\Z}$
of $\bar B'\setminus\sigma(B)$ and by gluing for all $i\in\Z$ the pseudoconvex component of 
the border $\partial A_{i-1}$ of $A_{i-1}$ to the pseudoconcave component of $\partial A_i$
by means of $\sigma\circ\pi$ again.
If one glues in this way only the copies $(A_i)_{i\le0}$
and then caps the pseudoconcave end of $A_0$ by a copy of $B$ using $\sigma$,
one obtains a non-compact complex surface $\hat X$, a holomorphic map
$p:\tilde X\to\hat X$, a $(-1)$-curve $\hat C$ on 
  $\hat X$ and a point $\hat O$ on $\hat C$
such that $p$ extends the identity map of $\cup_{i\le0}A_i$, $\hat C$
is the isomorphic image of a rational curve $C$ of $\tilde X$ through $p$
and $p$ restricts to an isomorphism $\tilde X\setminus p^{-1}(\hat O)\to \hat X\setminus\{\hat O\}$.
In fact $p^{-1}(\hat O)$ is the union of the infinitely many rational curves
appearing after $C$ on $\tilde X$ in the ``order of creation'', cf. \cite{Dl84}, Prop. 3.4.
Thus $p:\tilde X\to\hat X$ can be seen as a blowing down of the infinitely many exceptional curves
in $p^{-1}(\hat O)$. The generator $g$ of $\pi_1(X)$ mapping $A_i$ to $ A_{i+1}$ induces a holomorphic map
$\hat g:\hat X\to\hat X$ with $p\circ g= \hat g\circ p$.
One sees that $\hat O$ is fixed by $\hat g$, that $\hat g(\hat C)=\hat O$, and
 that the germ of $\hat g$ at $\hat O$
is the same as the germ of $\pi\circ\sigma:B\to B$ at the origin.

Let now $u$ be any plurisubharmonic function on $\tilde X$. 
The restriction of $u$ to $\hat X\setminus\{\hat O\}$
extends to a plurisubharmonic function $\hat u$ on $\hat X$, \cite{GR56}. 
It is clear that $\hat u\circ\hat g=\hat u -c$ in case $ u\circ g= u -c$ 
for some $c\in\R$. Conversely one gets a plurisubharmonic function
$u=\hat u\circ p$ starting from $\hat u$. In fact, since the germ of
$\hat g$ around $\hat O$ is contracting, it is enough to have only a germ of
$\hat u$ around $\hat O$ satisfying $\hat u\circ\hat g=\hat u -c$
in order to recover $u$ on $\tilde X$ with the property
 $ u\circ g= u -c$.

According to \cite{Dl84}, \cite{Dl88}, \cite{DO99}, \cite{Fav00} we get the following three
normal forms for representatives of conjugacy classes of Dloussky germs:
\begin{enumerate}
\item in the case of  Enoki surfaces
$$f(z,w)=(\alpha z,wz^s+Q(z)),$$ 
where $\alpha\in\Delta^*=\Delta\setminus\{0\}$, $s\ge1$ and $Q$ is a complex polynomial of degree at most
$s$ and with $Q(0)=0$; we have denoted by $\Delta$ the unit disc in $\C$; 
\item in the case of intermediate surfaces
$$f(z,w)=( z^p,\lambda wz^s+Q(z)),$$
 where $p\ge2$, $s\ge1$, $\lambda\in\C^*$ and $Q(z)=\sum_{m=1}^s a_mz^m+az^{\frac{ps}{p-1}}$
is a complex polynomial with $\gcd\{p, \ m \ | \ a_m\ne0\}=1$ and $a=0$ if $(p-1)\nmid s$ or $\lambda\ne 1$; 
\item in the case of Inoue-Hirzebruch surfaces
$$f(z,w)=(z^aw^b,z^cw^d),$$
where the matrix  $\begin{pmatrix} a& b\\c&d\end{pmatrix}$ is a product of 
$b_2(X)$ matrices of the form $\begin{pmatrix} 0& 1\\1&1\end{pmatrix}$
or $\begin{pmatrix} 1& 1\\0&1\end{pmatrix}$ with at least one factor of the first kind.
\end{enumerate}

In the rest of this paper we shall determine  the germs of plurisubharmonic functions
around the origin of $\C^2$ satisfying $ u\circ f= u -c$ for some fixed $c\in\R_{>0}$ 
and each type of germ $f$ as above.


\section{The main results}

\begin{Th}\label{enoki}
 (a) Let $f:\C^2\to\C^2$, 
$$f(z,w)=(\alpha z,wz^s+Q(z)),$$ 
with $\alpha\in\Delta^*$, 
$s\ge1$ and $Q$ is a complex polynomial of degree at most
$s$ and with $Q(0)=0$. Let $u:\C^2\to[-\infty,\infty[$, $u(z,w)=\log|z|$.
Then up to some additive constant, $u$ is the only plurisubharmonic function on
$\C^2$ which satisfies $u\circ f=u+\log|\alpha|$.

 (b) On an Enoki surface the integration current on the cycle of rational curves 
is the only positive exact current of bidimension $(1,1)$ 
up to multiplicative constants.
\end{Th} 
\begin{proof}
We start by proving part (b) of the theorem.

Let $C=C_1+...+C_n$ be the cycle of rational curves on the Enoki surface $X$. Since $C$ is homologically trivial the 
 current of integration $[C]$ along $C$ is a positive $d$-exact current of bidimension $(1,1)$.
We denote by $E$ the elliptic curve on $X$ in case it exists.

Let  $T$ be an arbitrary positive exact current of bidimension $(1,1)$ on $X$ and
$m_{C_i}:=\inf\{\nu(T,x) \ | x\in C_i\}$, $m_C=\min m_{C_i}$, $m_E$ its generic Lelong numbers along $C_i$, $C$ and $E$ respectively. 
We denote by $\chi_A$ the characteristic function of a subset $A$ of $X$.

The Siu decomposition of $T$ has the form
$$T= \chi_ET +\chi_CT+\chi_{X\setminus (C\cup E)}T
=m_E [E]+\sum_i m_{C_i}[C_i]+\chi_{X\setminus (C\cup E)}T=$$
$$=m_E [E]+
\sum_i ( m_{C_i}-m_C)[C_i]+m_C[C]+\chi_{X\setminus (C\cup E)}T,$$

see \cite{D} 6.18, 3.2.4. Since $C$ is the only homologically trivial effective reduced non-trivial divisor 
$X$ we get as in Remark \ref{rem} that $\chi_{X\setminus (C\cup E)}T$
is exact and that
$$T= m_C[C]+\chi_{X\setminus (C\cup E)}T.$$

Hence we may replace 
$T$ by $\chi_{X\setminus (C\cup E)}T$ which is positive, $d$-exact, has vanishing generic Lelong number
along $C$  and is the trivial extension to 
 $X$ of the restriction of $T$ to $X\setminus (C\cup E)$. After normalisation we obtain a
corresponding plurisubharmonic function $v$ on $\C^2$ satisfying $v\circ f=v+\log|\alpha|$.
We denote the current $dd^cv$ on $\C^2$ by $T$ again. Consequently we have
$f^*T=T$. Moreover, since $C$ corresponds to the axis $A=\{ z=0\}\subset\C^2$, the generic Lelong
number
$m_A$ of $T$ along $A$ must vanish. 
Since $f(A)=\{0\}$ and in general $\nu(f^*T,x)\ge\nu(T,f(x))$ (cf. \cite{Meo96}) it follows that
$\nu(T,0)=0$.

The differential form $dz$ defines a holomorphic foliation on   $\C^2$ which is invariant by $f$, since
$f^*(dz)=\alpha dz$.

{\em Claim}: $T$ is a foliation current for this foliation.

This means that for any test function $\phi\in \cal{C}^{\infty}_c(\C^2)$ one has
$T(\phi dz\wedge d\bar z)=0$. As before $T$ is the trivial extension to   $\C^2$  of its restriction to
 $\C^2\setminus A$. It is therefore enough to check that $ T(\phi dz\wedge d\bar z)=0$ for test
functions $\phi$ with support in $\C^2\setminus A$.

We have
$$|T(\phi dz\wedge d\bar z)|=|f^*T(\phi dz\wedge d\bar z)|=
|T((f^{-1})^*(\phi dz\wedge d\bar z))|=$$
$$=|T(\frac{\phi\circ f ^{-1}}{|\alpha|^2}dz\wedge d\bar z)|\le
\frac{\max |\phi|}{|\alpha|^2}2\sigma_T(f(\Supp \phi)),$$
where $\sigma_T=T\wedge\frac{i}{2}(dz\wedge d\bar z+dw\wedge d\bar w)$ denotes
the trace measure of $T$.
Iterating we obtain
$$|T(\phi dz\wedge d\bar z)|\le \frac{\max |\phi|}{|\alpha|^{2n}}2\sigma_T(f^n(\Supp \phi)), \leqno{(1)}$$
for any $n\in\N$.

We now need to estimate how large $f^n(B(0,R))$ is.
Take $C_1=\max\{\frac{|Q(z)|}{|z|} \ | \ z\in\bar\Delta\}$. We denote by $P(R_1,R_2)$ the bidisc of
radii $R_1$, $R_2$ centered at the origin of $\C^2$.
Then 

$$f(P(1,R_2)\subset P(|\alpha|, R_2+C_1),$$
$$f^2(P(1,R_2)\subset P(|\alpha|^2, |\alpha|(R_2+C_1)+|\alpha|C_1)\subset$$
$$\subset P(|\alpha|^2, |\alpha|(R_2+C_1)+(\frac{1}{1-|\alpha|}-1)C_1),$$
$$f^3(P(1,R_2)\subset P(|\alpha|^3, |\alpha|^3(R_2+C_1)+\frac{|\alpha|^2}{1-|\alpha|}|C_1)\subset$$
$$\subset P(|\alpha|^3, |\alpha|^3(R_2+C_1)+(\frac{1}{1-|\alpha|}-1)C_1),$$
$$f^4(P(1,R_2)\subset P(|\alpha|^4, |\alpha|^6(R_2+C_1)+\frac{|\alpha|^3}{1-|\alpha|}|C_1)\subset$$
$$\subset P(|\alpha|^4, |\alpha|^6(R_2+C_1)+(\frac{1}{1-|\alpha|}-1)C_1),$$
and further
$$f^n(P(1,R_2)\subset P(|\alpha|^n, |\alpha|^{\frac{n(n-1)}{2}}(R_2+C_1)+\frac{|\alpha|^{n-1}}{1-|\alpha|}|C_1).$$

From this and from (1) it follows that there is a constant $C_2$ such that 
$$|T(\phi dz\wedge d\bar z)|\le C_2\max |\phi|\frac{\sigma_T(B(0,|\alpha|^n )}{|\alpha|^{2n}},$$
for all $n$ sufficiently large.
But $$\lim_{n\to\infty}\frac{\sigma_T(B(0,|\alpha|^n )}{|\alpha|^{2n}}=\pi\nu(T,0)=0,$$
and the claim follows.

Using the claim it can be easily shown that $T$ is invariant under translations in the
$w$-direction, hence $T=pr_1^*S$, where $S$ is a positive current on $\C$ of dimension zero and
$pr_1:\C^2\to\C$ denotes the first projection. Set $f_1:\C\to\C$, $f_1(z)=\alpha z$.
Then $pr_1\circ f=f_1\circ pr_1$, hence $f_1^*(S)=S$. The current $S$ is of the form 
$\mu idz\wedge d\bar z$, where $\mu$ is a positive measure on $\C$. 
Denote by $\Delta(r)$ the disc of radius $r$ around the origin of $\C$. 
The invariance property
of $S$ implies
$$\mu(\Delta(r))=\mu(\Delta(|\alpha|r))=\mu(\Delta(|\alpha|^nr))$$
for all $n\in\N$ and $r\in\R_{>0}$. This entails that $S$ is supported at the origin
of $\C$ and hence that $T=0$, since $T$ is not carried by $A$.\\

We now turn to part (a) of the theorem.
Let $v$ be a plurisubharmonic function on $\C^2$ satisfying 
$v\circ f=f+\log|\alpha|$. By part (b) of the theorem we see that
$dd^cv\le dd^cu$ or $dd^cu\le dd^cv$. But then $u-v$ or $v-u$ would give a plurisubharmonic function on 
$X$ which has to be constant. 
\end{proof}

{\em Notation} For $\alpha\in\R_{>0}$ we denote by $K_{\alpha}$ the set of continuous 
$\alpha$-periodic functions $\psi:\R\to\R$ which fulfil the inequality
$$-\psi''+\psi'+1\ge0$$
in generalised sense.

Notice that $K_{\alpha}$ is infinite dimensional: for any $\alpha$-periodic smooth function $\phi:\R\to\R$ and for any small enough factor $\epsilon\in\R_{>0}$ one has 
$\epsilon\phi\in K_{\alpha}$.

\begin{Th}\label{intermed}
 (a) Let $f:\Delta\times\C\to\Delta\times\C$, 
$$f(z,w)=( z^p,\lambda wz^s+Q(z)),$$
 where $p\ge2$, $s\ge1$, $\lambda\in\C^*$ and $Q(z)=\sum_{m=1}^s a_mz^m+az^{\frac{ps}{p-1}}$
is a complex polynomial with $\gcd\{p, \ m \ | \ a_m\ne0\}=1$ and $a=0$ if $(p-1)\nmid s$ or $\lambda\ne 1$.
The  plurisubharmonic functions
 $u:\Delta\times\C\to[-\infty,\infty[$
 which satisfy $u\circ f=u-\log p$ are precisely the functions of the form
$$u(z,w)=-\log(-\log|z|)-\psi(\log(-\log|z|))$$
for $\psi\in K_{\log p}$.

 (b) The cone of positive exact currents of bidimension $(1,1)$ on the intermediate surface
$X(f)$ corresponds bijectively to the cone of currents of the form $cdd^c u$ on $\Delta\times\C$
for $u$ as above and $c\in\R_{>0}$. In particular,
the modified K\"ahler rank of $X(f)$ is one.
\end{Th} 
\begin{proof}
There are no
 homologically trivial divisors on $X(f)$ this time (excepting $0$ of course). 
As in the case of Enoki surfaces we reduce ourselves to the
investigation of closed positive currents 
  $T$ on $\Delta\times\C$ with vanishing Lelong number
 at the origin and satisfying $f^*T=T$. Moreover we may again suppose that  $T$ is the extension to  $\Delta\times\C$ 
of its restriction to $\Delta^*\times\C$.

The differential form $dz$ defines a holomorphic foliation on $\Delta\times\C$ which is invariant 
under $f$. 

We start again by showing that $T$ is a
foliation current for this foliation. For this it is enough to
check that the measure

$$idz\wedge d\bar z\wedge T$$
vanishes on $\Delta^*\times\C$. We use the invariance of $T$ by $f$ again.

Take $0<r_1<r_2<1$, $r'>0$, $D:=(\Delta(r_2)\setminus\Delta(r_1))\times \Delta(r')$, 
$D':=f^{-n}(f^n(D))$ and 
$A(r_1,r_2)=A(r_1,r_2,r'):=(i
dz\wedge d\bar z\wedge T)(D)$
the measure of the set $D$. We have of course $D'\supset D$. Since $f:\Delta^*\times\C\to\Delta^*\times\C$
is $p$ to $1$ we get
$$(i
dz\wedge d\bar z\wedge T)(f^n(D))=
p^n(i|z|^{2(p^{n}-1)}dz\wedge d\bar z\wedge T)(D')\ge$$
$$\ge p^n(i|z|^{2p^{n}}
dz
\wedge d\bar z\wedge T)(D')\ge
p^n r_1^{2p^{n}}A(r_1,r_2).$$
 Hence

$$A(r_1,r_2)\le \frac{1}{p^n}(\frac{r_2}{r_1})^{2p^{n}}
\frac{i(\partial\bar\partial(|z|^2+|w|^{2p})\wedge T)(f^n(D))}{\pi r_2^{2p^{n}}}. \leqno{(2)}$$

As before we need to estimate the width of $f^n(D)$. 
Let $C_1=\max\{\frac{|Q(z)|}{|z|} \ | \ z\in\bar\Delta\}$ and $C_2=\max(1,|\lambda|).$
For $0<r<1$, $0<r'$ we see that
$$f(P(r,r'))\subset P(r^p,rC_2(r'+\frac{C_1}{C_2}))\subset  P(r^p,rC_2(r'+\frac{C_1}{(1-r)C_2})),$$
$$f^2(P(r,r'))\subset P(r^{p^2},r^{p+1}C_2^2r'+\frac{r^{p-1}C_2^2C_1}{(1-r)C_2})\subset  P(r^{p^2},r^{p-1}C_2^2(r'+\frac{C_1}{(1-r)C_2}),$$
$$f^3(P(r,r'))\subset P(r^{p^3},r^{p^2}C_2(r^{p-1}C_2^2(r'+\frac{C_1}{(1-r)C_2})+\frac{C_1}{C_2}))\subset $$ $$\subset P(r^{p^3},r^{p^2+p-1}C_2^3r'+\frac{r^{p^2-1}C_2^3 C_1}{(1-r)C_2}).$$
For sufficiently large $n$ we thus obtain

$$f^n(P(r,r'))\subset\{ |z|^2+|w|^{2p}<2C_1 r^{p^{n}}\}$$

which in combination with the inequality (2) gives

$$A(r_1,r_2)\le \frac{1}{p^n}(\frac{r_2}{r_1})^{2p^{n}}
\frac{2i(\partial\bar\partial(|z|^2+|w|^{2p})\wedge T)(\{ |z|^2+|w|^{2p}<2C_1 r_2^{p^{n}}\})}
{\pi r_2^{2p^{n}}}.$$

But 
 the factor 
 $$\frac{2i(\partial\bar\partial(|z|^2+|w|^{2p})\wedge T)(\{ |z|^2+|w|^{2p}<2C_1 r_2^{p^{n}}\})}
{\pi r_2^{2p^{n}}}$$
 converges to $4C_1 \nu(T,\phi,0)$, where  $\nu(T,\phi,0)$ is the Lelong number
 with respect to $\phi=\log(|z|^2+|w|^{2p})$. This Lelong number vanishes
 by the Comparison Theorem for Lelong numbers \cite{D}, since the usual Lelong number vanishes. 
 So for 
any $\epsilon>0$ we can find some $N\in\N$ such that

$$A(r_1,r_2)\le \frac{1}{p^n}(\frac{r_2}{r_1})^{2p^{n}}\epsilon$$ 
for all $n\ge N$. Moreover this inequality holds also for smaller  $r_1$ and $r_2$.

Let 
$$\delta=(\frac{r_2}{r_1})^{\frac{1}{2p^n}}$$
and look at the division $(r_1,\delta r_1, \delta^2 r_1,...,\delta^{2p^n}r_1=r_2) $
of the interval $[r_1,r_2]$.
We have seen that 
$A(\delta^ir_1,\delta^{i+1}r_1)\le \frac{1}{p^n} \delta^{2p^n}\epsilon=\frac{1}{p^n}\frac{r_2}{r_1}\epsilon$
hence 
$A(r_1,r_2)=\sum_{i=1}^{2p^n}A(\delta^{i-1}r_1,\delta^{i}r_1)\le 2\frac{r_2}{r_1}\epsilon$,
which proves that $T$ is a foliation current.

As in the proof of Theorem \ref{enoki}, $T$ is invariant under translations in the
$w$-direction, hence $T=pr_1^*S$, where $S$ is a positive current on $\Delta$ of dimension zero and
$pr_1:\Delta\times\C\to\Delta$ is the first projection. On the other hand $T=dd^cu$, where u is a
plurisubharmonic function on $\Delta\times\C$ with $u\circ f= u-\log p$.
The restriction of $u$ to a general leaf $\{ z=const.\}$ of our foliation is a harmonic function in
$w$. We shall show that $\frac{\partial u}{\partial w}=0$, the case of 
$\frac{\partial u}{\partial\bar w}=0$ being similar.

Since $T$ is a foliation current one has that $\frac{\partial^2 u}{\partial\bar z\partial w}=
\frac{\partial^2 u}{\partial z\partial w}=0$ as distributions. In particular $\frac{\partial u}{\partial w}=0$
is invariant by translations in the $z$-directions.
The invariance condition by translations in $w$ implies that
$$\lambda z^s\frac{\partial u}{\partial w}(f(z,w))=\frac{\partial u}{\partial w}(z,w). \leqno{(3)}$$
Take now a small non-zero solution $z'$ of the equation
$$ w=\lambda z^sw+Q(z) \leqno{(4)}$$
for fixed small $w$, $w\ne 0$. Using the invariance in the $z$-direction together with (3) and (4)
we obtain

$$\frac{\partial u}{\partial w}(z,w)=\frac{\partial u}{\partial w}(z',w)
=\lambda (z')^s\frac{\partial u}{\partial w}((z')^p,w)=\lambda (z')^s\frac{\partial u}{\partial w}(z,w),$$
hence $\frac{\partial u}{\partial w}(z,w)=0$ for small $w$. But for large $w$ the vanishing holds as well since we
can iterate on (3).

Thus $u$ is constant on the leaves $\{ z=const.\}$ and therefore descends to a subharmonic function
$u$ on $\Delta$ fulfilling the condition
$u(z^p)=u(z)-\log p$ for all $z\in\Delta$.

Take now $z\in\Delta$ and a $p^n$-th root of unity $\theta$. The invariance condition
implies $u(z)=u(z^{p^n})+n\log p=u(\theta z)$ hence each value $u(z)$ is attained on a dense subset of the circle
$\{ |z|=const.\}$. By semi-continuity it follows that $u$ depends only on $r=|z|$.

Let $v=u|_{[0,1[}$. By the maximum principle and the invariance condition on $u$
one sees that $v$ is strictly increasing. Using this, the semi-continuity of $u$ 
and the definition of subharmonicity (cf. \cite{D} 1.4.13.b) one infers that $v$
is continuous.
The Laplace operator takes the form 
$\frac{\partial^2}{\partial r^2}+\frac{1}{r}\frac{\partial}{\partial r}
+\frac{1}{r^2}\frac{\partial^2}{\partial \phi^2}$ in polar coordinates $(r,\phi)$.
We thus reduce ourselves to the search of continuous functions
$v:[0,1[\to[-\infty,\infty[$ which satisfy

$$v(r^p)=v(r)-\log p \leqno{(5)}$$

and

$$v''(r)+\frac{1}{r}v'(r)\ge0. \leqno{(6)}$$

Let $h:]-\infty,\infty]\to[0,1[$, $h(t)=\exp(-\exp t)$
and $\psi=-v\circ h-h$.
The conditions (5) and (6) translate into

$$\psi(t+\log p)=\psi(t),$$
$$-\psi''+\psi'+1\ge0$$

and thus $\psi\in K_{\log p}$ and the theorem is proved.
\end{proof}

For the case of Inoue-Hirzebruch surfaces we need some preparations.
Let $f:\C^2\to\C^2$, $f(z,w)=(z^aw^b,z^cw^d)$, and denote by $A$
 the matrix  $\begin{pmatrix} a& b\\c&d\end{pmatrix}$. We suppose that
$A$  is a product of 
matrices of the form $\begin{pmatrix} 0& 1\\1&1\end{pmatrix}$
or $\begin{pmatrix} 1& 1\\0&1\end{pmatrix}$ with at least one factor of the first kind.
We have $\det(A)=\pm 1$ and by Lemma 2.2 from \cite{Dl88}, $\trace(A)>2$ unless
$A=\begin{pmatrix} 0& 1\\1&1\end{pmatrix}$
or $A=\begin{pmatrix} 0& 1\\1&2\end{pmatrix}$.
In any case one sees that the eigenvalues $\lambda_1$, $\lambda_2$ of $A^t$ are real, 
irrational and one of them is larger than $1$. 
We set $\lambda=\lambda_1>1$. 
Let $(\alpha,\beta)=(\alpha_1,\beta_1)$, $(\alpha_2,\beta_2)$ be eigenvectors
associated to $\lambda_1$ and $\lambda_2$ respectively.
An easy computation shows that $\alpha\beta=\alpha_1\beta_1>0$, $\alpha_2\beta_2<0$.
We choose $\alpha,\beta,\alpha_2\in\R_{>0}$ and set
$$\phi=\phi_1=\alpha\log|z|+\beta\log|w|, \ \phi_2=\alpha_2\log|z|+\beta_2\log|w|,$$
$$U:=\{\phi<0\}\subset\C^2.$$
Then $\phi_i\circ f=\lambda_i\phi_i$ for $i=1,2$.

Now we can state
\begin{Th}\label{IH}
 (a) For $A$, $\lambda$, $\alpha$, $\beta$ as above let
$\phi=\alpha\log|z|+\beta\log|w|$, $U=\phi^{-1}([-\infty,0[)$,
 $f:U\to U$, 
$$f(z,w)=(z^aw^b,z^cw^d).$$
Then the  plurisubharmonic functions
 $u$ on $U$
 which satisfy $u\circ f=u-\log \lambda$ are precisely the functions of the form
$$u=-\log(-\phi)-\psi(\log(-\phi)),$$
for $\psi\in K_{\log p}$.

 (b) The cone of positive exact currents of bidimension $(1,1)$ on the Inoue-Hirzebruch surface
$X(f)$ corresponds bijectively to the cone of currents of the form $cdd^c u$ on $U$
for $u$ as above and $c\in\R_{>0}$. In particular,
the modified K\"ahler rank of $X(f)$ is one.
\end{Th} 
\begin{proof}

As in the proof of Theorem \ref{intermed} we reduce ourselves to the
investigation of currents 
  $T$ on $U$ of the form $T=dd^c u$ with $u$ plurisubharmonic on $U$ and such that
$u\circ f=u-\log \lambda$
and $\nu(T,0)=0$. Moreover we may  suppose that  $T$ is the extension to  $U$ 
of its restriction to $U\cap(\C^*\times\C^*)$.

There is again a holomorphic foliation on $U$, this time singular, which will be shown
 to be compatible with $T$.
On $U\cap(\C^*\times\C^*)$ this foliation is induced by the smooth positive $(1,1)$-form 
$dd^cv$ where $v:=-\log(-\phi)$. As before it is enough to show that the measure 
$dd^cv\wedge T$ is zero on any compact subset of $U\cap(\C^*\times\C^*)$.

For a better visualisation of the present situation we introduce the map
$E:\C\times\C\to\C^*\times\C^*$, $E(\zeta,\omega)=(\exp\zeta,\exp\omega)$.
We denote the linear automorphism induced by $A$ on $\C^2$ again by $A$ and
get the relations
$$f\circ E=E\circ A,$$
$$\phi\circ E(\zeta,\omega)=\alpha_i\Re e\zeta+\beta_i\Re e\omega, \ i=1,2.$$
In particular $E^{-1}(U)=\{\alpha\Re e\zeta+\beta\Re e\omega<0\}$.

We also see that $U$ is covered by compact subsets of the form
$D=D(c_1,\delta,c_2):=\phi_1^{-1}([-c_1\delta,-c_1])\cap\phi_2^{-1}([-c_2,c_2])$ for
$c_1,c_2\in\R_{>0}$, $\delta\in[1,\infty[$.
Since both $T$ and $dd^cv$ are  invariant by $f$ we have
$$A(c_1,\delta)=A(c_1,\delta,c_2):=(dd^cv\wedge T)(D(c_1,\delta,c_2))=$$
$$=(dd^cv\wedge T)f^n(D(c_1,\delta,c_2)))=(dd^cv\wedge T)(D(\lambda^nc_1,\delta,\lambda^{-n}c_2)).$$

We need to estimate $dd^cv$ on $f^n(D)$. On $D$ we have

$$-c_1\delta\le\alpha_1\log|z|+\beta_1\log|w|\le-c_1,$$
$$-c_2\le\alpha_2\log|z|+\beta_2\log|w|\le c_2,$$

hence

$$\exp(\frac{c_1\delta\beta_2-c_2\beta_1}{\alpha_2\beta_1-\alpha_1\beta_2})
\le|z|\le
\exp(\frac{c_1\beta_2+c_2\beta_1}{\alpha_2\beta_1-\alpha_1\beta_2}),$$
$$\exp(\frac{-c_1\delta\alpha_2-c_2\alpha_1}{\alpha_2\beta_1-\alpha_1\beta_2})
\le|w|\le
\exp(\frac{-c_1\alpha_2+c_2\alpha_1}{\alpha_2\beta_1-\alpha_1\beta_2}).$$

Therefore we get on $f^n(D)$

$$\exp(\frac{\lambda^nc_1\delta\beta_2-\lambda^{-n}c_2\beta_1}{\alpha_2\beta_1-\alpha_1\beta_2})
\le|z|\le
\exp(\frac{\lambda^nc_1\beta_2+\lambda^{-n}c_2\beta_1}{\alpha_2\beta_1-\alpha_1\beta_2}),$$
$$\exp(\frac{-\lambda^nc_1\delta\alpha_2-\lambda^{-n}c_2\alpha_1}{\alpha_2\beta_1-\alpha_1\beta_2})
\le|w|\le
\exp(\frac{-\lambda^nc_1\alpha_2+\lambda^{-n}c_2\alpha_1}{\alpha_2\beta_1-\alpha_1\beta_2}),$$

hence $f^n(D)\subset B(0,r_n)$, where

$$r_n^2=\exp(2\frac{\lambda^nc_1\beta_2+\lambda^{-n}c_2\beta_1}{\alpha_2\beta_1-\alpha_1\beta_2})+
\exp(2\frac{-\lambda^nc_1\alpha_2+\lambda^{-n}c_2\alpha_1}{\alpha_2\beta_1-\alpha_1\beta_2}).$$

Now on $D$ again
$$i\partial\bar\partial v=\frac{i\partial\phi\wedge\bar\partial\phi}{\phi^2}=
\frac{i(\alpha\frac{dz}{z}+\beta\frac{dw}{w})\wedge(\alpha\frac{d\bar z}{\bar z}+\beta\frac{d\bar w}{\bar w})}
{4\phi^2}\le$$
$$\le\frac{C_1}{\min\{|z|^2,|w|^2 \ | \ (z,w)\in D\}}\frac{i(dz\wedge d\bar z+dw\wedge d\bar w)}{\phi^2},$$

where $C_1$ is a constant not depending on $D$.
This implies 

$$A(c_1,\delta)\le
\frac{C_1r_n^2}{\min\{\phi^2(z,w)|z|^2,\phi^2(z,w)|w|^2 \ | \ (z,w)\in f^n(D)\} }
\frac{\sigma_T(B(0,r_n^2))}{r_n^2}.$$

As in the previous cases $\frac{\sigma_T(B(0,r_n^2))}{r_n^2}$ converges to zero, therefore for any
$\epsilon>0$ we get 
$\frac{C_1\sigma_T(B(0,r_n^2))}{r_n^2}<\epsilon$ as soon as $n$ is large enough.
Setting $C_2=\frac{2c_1\max\{\alpha_1,-\beta_2\}}{\alpha_2\beta_1-\alpha_1\beta_2}+1$ we get
by our estimates
$r_n^2\le\exp(-\lambda^nC_2)$ and
$\min\{|z|^2,|w|^2 \ | \ (z,w)\in f^n(D)\}\ge\exp(-\lambda^nC_2\delta)$,
hence

$$A(c_1,\delta)\le\frac{\exp(\lambda^nC_2(\delta-1))}{\lambda^{2n}}\frac{\epsilon}{c_1^2}.$$

Set $m=\lfloor\lambda^{2n}\rfloor$ and consider the division 
$(1, \exp\frac{1}{m}, \exp\frac{2}{m},...,\exp\frac{m}{m}=e)$ of the interval 
$[1,e]$.
For any  subinterval our estimates hold with $eC_2$ instead of $C_2$, hence

$$A(c_1,e,c_2)=\sum_{i=1}^mA(c_1\exp\frac{i-1}{m},\exp\frac{1}{m},c_2)\le$$
$$\le m\frac{\exp(\lambda^nC_2(\exp\frac{1}{m}-1))}{\lambda^{2n}}\frac{\epsilon}{c_1^2}
\le \exp(C_2e\frac{\exp\frac{1}{m}-1}{\frac{1}{m}})\frac{\epsilon}{c_1^2}$$

and the last term converges to $\exp(C_2e)\frac{\epsilon}{c_1^2}$ showing that
$A(c_1,e,c_2)=0$.

As in the proof of Theorem \ref{intermed} we next check that $u$ is constant
on the leaves of the foliation given by $\alpha\frac{dz}{z}+\beta\frac{dw}{w}$
on $U\cap(\C^*\times\C^*)$. 

Indeed, for $V:=u\circ E$ one has $V\circ A=V-\log \lambda$ on
$\{\alpha\Re e\zeta+\beta\Re e\omega<0\}$ and $dd^cV$ is a foliation current for the foliation
given by $\alpha d\zeta+\beta d\omega$.

We consider a linear change of coordinates which diagonalizes $A$ on $\C^2$ leading to
$A(\xi,\tau)=(\lambda\xi,\lambda^{-1}\tau)$. Keeping the notation $V$ for $V$ after this coordinate
change we notice that $V$ restricted to the leaves $\{\xi=const.\}$
 is harmonic and that $\frac{\partial V}{\partial \tau}$, $\frac{\partial V}{\partial \bar\tau}$
are invariant by translations in the $\xi$-direction.
As in the case of intermediate surfaces we get
$$\frac{1}{\lambda}\frac{\partial V}{\partial \tau}(\xi,\frac{1}{\lambda}\tau)=
\frac{1}{\lambda}\frac{\partial V}{\partial \tau}(\lambda\xi,\frac{1}{\lambda}\tau)=
\frac{\partial V}{\partial \tau}(\xi,\tau).$$ 
Iterating this relation and using the continuity of 
$\frac{\partial V}{\partial \tau}$ on $\{\xi=const.\}$ we obtain
$$\frac{\partial V}{\partial \tau}(\xi,\frac{1}{\lambda^n}\tau)=\lambda^n\frac{\partial V}{\partial \tau}(\xi,\tau)$$
which forces $\frac{\partial V}{\partial \tau}$ to vanish.
Similarly $\frac{\partial V}{\partial \bar\tau}=0$. Thus $V$ is constant on the leaves $\{\xi=const.\}$.

Returning now to the coordinates $(\zeta,\omega)$ on $\C^2$ we shall show that $V$ only depends on 
$\Re e\zeta$ and $\Re e\omega$.
The function $V$ is doubly periodic in $\Im m\zeta$, $\Im m\omega$ since $V=u\circ E$.
Fix $\Re e\zeta$ and $\Re e\omega$ and look at a leaf $\{\alpha \zeta+\beta \omega=const.\}$.
Under this restriction the imaginary parts must satisfy some relation 
$\alpha\Im m \zeta+\beta\Im m \omega=const.$ which describes a dense subset in the ``torus of the
imaginary parts'' $\R^2/(2\pi\Z)^2$ since $\alpha/\beta$ is irrational.
On this subset $V$ is constant and by semi-continuity it must be constant on the whole torus.

Switching once more to coordinates 
$(\xi,\tau)$ we see that $V$ depends on $\Re e\xi$ alone.
We have reduced ourselves in this way to the search of  continuous functions 
$v:\R_{<0}\to\R$ subject to the conditions
$$v(\lambda t)=v(t)-\log \lambda, \ \forall t\in\R_{<0},$$
$$v''\ge0.$$

Remark that the continuity of $v$ can be deduced as in Theorem \ref{intermed} after first
restricting $u$ to the line $\{ z=w \}$.

Take now $h:\R\to\R_{<0}$, $h(s)=-\exp(s)$ and $\psi=-v\circ h-h$. Then the conditions on $v$ translate into
$\psi\in K_{\log\lambda}.$
Thus $u=-\log(-\phi)-\psi(\log(-\phi))$ and the proof is finished.
\end{proof}

Finally, the first part of Theorem \ref{foliated} is a direct 
consequence of the description of $P_{bdy}(X)$ given in \cite{HL83} and in the above Theorems.
For the second part it suffices to apply Remark \ref{rem} and to notice that the only case of a curve on a ``known'' non-K\"ahlerian surface, which is not invariant under a holomorphic foliation, is that of the elliptic curve on a parabolic Inoue surface.


\begin{thebibliography}{123}


\bibitem{BHPV} {\sc Barth W., Hulek K., Peters C., Van de Ven A.}:
Compact complex surfaces, {\em Ergebnisse der Mathematik und ihrer Grenzgebiete, 
Berlin, Springer-Verlag,  (2004).}

\bibitem{D} {\sc Demailly J.P.}: Complex analytic and algebraic geometry,
http://www-fourier.ujf-grenoble.fr/$\sim$demailly/books.html.

\bibitem{Dl84} {\sc Dloussky G.}: Structure des surfaces de
Kato, {\em M\'emoire de  la S.M.F 112.$\rm n^{\circ}14$ (1984).}  
\bibitem{Dl88} {\sc Dloussky G.}: Sur la classification des germes d'applications holomorphes contractantes,
 {\em Math. Ann. 280, (1988), 649-661.} 
\bibitem{DO99} {\sc Dloussky G., Oeljeklaus K.}: Vector fields and foliations on surfaces of class VII${}_{0}$,
 {\em Ann. Inst. Fourier 49, (1999), 1503-1545}.
\bibitem{Fav00} {\sc Favre, Ch.}: Classification of $2$-dimensional
contracting rigid germs, 
{\em Jour. Math. Pures Appl. 79 (2000), 475-514.}
\bibitem{GR56} {\sc Grauert, H., Remmert, R.}: Plurisubharmonische Funktionen in komplexen R\"aumen,
{\em Math. Z. 65 (1956), 175-194.}
\bibitem{HL83}{\sc Harvey, R., Lawson, H. B.}:
An intrinsic characterization of K‰hler manifolds,
{\em Invent. Math. 74 (1983), 169-198.}
\bibitem{Kato78} {\sc Kato Ma.}: Compact complex manifolds containing ``global spherical
shells'', {\em Proceedings of the Int. Symp. Alg. Geometry,  Kyoto 1977, Kinokuniya Book Store, Tokyo 1978.}
\bibitem{Lam99} {\sc Lamari, A.}: Courants k\"ahl\'eriens et surfaces compactes,
{\em  Ann. Inst. Fourier 49 (1999), 263-285.} 
\bibitem{Lam99'} {\sc Lamari, A.}: Le c\^one K\"ahl\'erien d'une surface,
{\em   J. Math. Pures Appl. 78 (1999), 249-263}.

\bibitem{Meo96} {\sc M\'eo, M.}:
   Image inverse d'un courant positif ferm\'e par une application analytique surjective, 
   {\em C. R. Acad. Sci. Paris SÈr. I Math. 322 (1996), 1141--1144.}
\bibitem{Nak90} {\sc Nakamura I.}: On surfaces of class $\rm VII_0$ with curves II,
{\em Toh\^oku Math. Jour. 42, (1990), 475-516.}

 \end{thebibliography}
\end{document}